\documentclass{article}
\usepackage{amsmath}
\usepackage{amscd}
\usepackage{latexsym}
\usepackage{amssymb}

\title{A note on the cohomology of Lie algebras}
\author{M. Gerstenhaber\footnote{The author wishes to thank V. Coll for helpful remarks during the preparation of this paper.}}

\begin{document}
\maketitle
\newtheorem{theorem}{Theorem}
\newtheorem{corollary}{Corollary}
\newtheorem{lemma}{Lemma}
%{\theorembodyfont{\rmfamily} \newtheorem{Rem}{Remark}}
\renewcommand{\abstractname}{}
\newcommand{\C}{\ensuremath{\mathbb{C}}}
\newcommand{\Z}{\ensuremath{\mathbb{Z}}}
\newcommand{\I}{\ensuremath{\mathcal{P}}}
\newcommand{\pr}{\ensuremath{\preceq}}
\newcommand{\op}{\ensuremath{\mathrm{op}}}
\newcommand{\g}{\mathfrak{g}}
\newcommand{\h}{\mathfrak{h}}
\newcommand{\G}{\mathfrak{G}}
\newcommand{\n}{\mathfrak{n}}
\newcommand{\m}{\mathfrak{m}}
\renewcommand{\b}{\mathfrak{b}}
\newcommand{\ad}{\ensuremath\operatorname{ad}}
\newcommand{\diag}{\ensuremath\operatorname{diag}}
\vspace{-2mm}
\date{}
{\noindent \textit{Department of Mathematics, University of
Pennsylvania, Philadelphia, PA 19104-6395}
\vspace{-7mm}

\begin{abstract}\noindent  This note presents a  general theorem about the cohomology of finite dimensional Lie algebras of arbitrary characteristic. As an application we compute the cohomology of the Borel subalgebra of $\mathfrak{sl}(N)$.
\end{abstract}

\section{A basic theorem} 
Suppose that $\g$ is a finite dimensional Lie algebra over a field $k$, which may be of any characteristic, and that $M$ is a finite dimensional $\g$ module. We will call an element $\sigma\in\g$ \emph{ad-semisimple (adss)} if it acts semisimply both on $\g$ and on $M$. (Note that unless $\g$ itself is semisimple, an element whose adjoint on $\g$ is semisimple need not act in a semisimple way on every finite dimensional representation.) Suppose now that $k$ contains all the eigenvalues of $\ad\sigma$, so that both $\g$ and $M$ decompose into a direct sum of eigenspaces.  If $\g_{\lambda}, \g_{\mu}, M_{\nu}$ are respective eigenspaces for eigenvalues $\lambda, \mu, \nu$ then $[\g_{\lambda}, \g_{\mu}] \subset \g_{\lambda+\mu}, [\g_{\lambda}, M_{\nu}] \subset M_{\lambda +\nu}$, so $\sigma$ induces a grading on both $\g$ and $M$ with values in the additive group of $k$. If either $M=\g$, or $M = k$ considered as a trivial $\g$ module, then $\sigma \in \g$  is already \emph{adss} if it is so on $\g$, in the first case by definition and in the second because the operation is identically zero for all elements of $\g$. Also, the \emph{adss} element $\sigma$ used to define the grading is always homogeneous of degree zero in the grading which it defines since $[\sigma, \sigma] = 0$, as are the elements of $k$ when $M=k$.

The grading induced by an \emph{adss} element extends to the cohomology of $\g$ with coefficients in $M$.  Recall that the latter may be defined, using the Chevalley-Eilenberg complex, as follows.
Suppose that $\g$ is a Lie algebra over an arbitrary commutative unital ring $k$ and that $M$ is a $\g$ module.  An $n$-cochain of $\g$ with coefficients in $M$ is a linear map $F^n: \bigwedge^n\g \to M$. Denote the $k$-module of these $n$-cochains by $C^n(\g, M)$ or simply $C^n$.  Then  the coboundary mapping $\delta: C^n \to C^{n+1}$ is given by 
\begin{multline*}
(\delta F)(a_0, a_1, \dots, a_n)=\\
\sum_{0\le i\le n}(-1)^i[ a_i, F(\dots,\hat a_i,\dots)] +
\sum_{0\le i < j \le n}(-1)^{i+j}F([a_i,a_j],\dots, \hat a_i,\dots, \hat a_j, \dots)
\end{multline*}
where $\hat a$ indicates omission of the argument $a$.  One has $\delta\delta =0$ so one defines, as usual, the $k$-modules of $n$-cocycles $Z^n = Z^n(\g,M)$, $n$-coboundaries $B^n = B^n(\g,M)$ and the $n$th cohomology group $H^n = H^n(\g,M) = Z^n/B^n$.  Here $C^0(\g, M)$ is just $M$ and  $H^0(\g, M) = Z^0(\g,M) =M^{\g}$, the set of all $m \in M$ such that $[a,m] = 0$ for all $a\in \g$, so $H^0(\g, k) = k$.

It will be useful to rearrange the terms on the right in the coboundary formula by taking first all those in which $\ad a_0$ appears as an operator on $M$ or on $\g$.  Denoting by  $\iota_{a_0}\!\cdot F$ the $n-1$-cochain defined by setting  $\iota_{a_0}\!\cdot F(a_1,\dots,a_{n-1} = F(a_{i_0}, a_1,\dots,a_{n-1})$, the coboundary then takes the form
\begin{multline}\label{coboundary}
(\delta F)(a_0, a_1, \dots, a_n)=\\
[a_0, F(a_1,\dots,a_n)] + \sum_{1 \le i \le n}(-1)^i F([a_0,a_i],a_1,\dots, \hat a_i, \dots, a_n)\\
-(\delta (\iota_{a_0}\!\!\cdot F))(a_1,\dots, a_n).
\end{multline}
 This permits a recursive definition with which it is easy to prove exactness, cf. e.g. \cite{HiltonStammbach}. 

Return now  to the case where $\g, \, M$ are finite dimensional over a field $k$ containing the eigenvalues of the adjoint of an \emph{adss} element $\sigma$. The space of $n$-cochains  $C^n(\g, M)$ then also acquires a grading: Define $F \in C^n$ to have degree $r$ if  $a_i \in \g_{\lambda_i}, \, i = 1,\dots,n$ implies $F(a_1, \dots, a_n) \in M_{\Sigma \lambda_i +r}$.  Every $n$-cochain is then the sum of its homogeneous parts.  Moreover, the grading induced by an \emph{adss} element is preserved by the coboundary operator. That is, if $F\in C^n$ is homogeneous of degree $r$ then so is $\delta F$.  The complex therefore decomposes into a direct sum of subcomplexes each consisting of the cochains of a fixed degree, so the cohomology  $H^*(\g, M)$ likewise decomposes into the sum of its homogeneous parts. This depends, of course, on the choice of $\sigma$.  Nevertheless, with these notations, one has the following

\begin{theorem}  Let $\g$ be a finite dimensional Lie algebra over a field $k$ and $M$ be a finite dimensional $\g$ module. Suppose that $\sigma$ is an \emph{adss} element of $\g$, that $k$ contains the eigenvalues of $\ad\sigma$, and that  $F$ is an $n$-cocycle in $C^n(\g,M)$ which is homogeneous of degree $r$ with respect to the grading induced by $\sigma$. 
\begin{enumerate}
\item If $r \ne 0$ then $F$ is a coboundary.
\item If $r = 0$ and $n \ge 1$ then the $n-1$ cochain $\iota_{\sigma}\!\cdot F$ is again  a 
 cocycle.
\end{enumerate}
\end{theorem}\label{zero}
\noindent\textsc{Proof.}  Suppose that $a_i, \, i = 1, \dots n$ are homogeneous elements of degrees $\lambda_i, i = 1,\dots, n$ respectively. Taking $a_0$ to be the \emph{adss} element $\sigma$, note that 
$$F([a_0,a_i], a_1,\dots,\hat a_i,\dots,a_n) = \lambda_iF(a_i,a_1,\dots,\hat a_i, \dots, a_n),$$
 so $(-1)^iF([a_0,a_i], a_1,\dots, \hat a_i,\dots, a_n)$ is just $-\lambda_iF(a_1,\dots,a_n)$. Therefore, $$\sum_{1 \le i \le n}(-1)^i F([a_0,a_i],a_1,\dots, \hat a_i, \dots, a_n) = -\sum_{1 \le i \le n} \lambda_i F(a_1,\dots a_n).$$  Since $F$ is homogeneous of degree $r$, its coboundary collapses to
$$(\delta F)(a_0, a_1,\dots a_n) = rF(a_1,\dots,a_n)-(\delta(\iota_{a_0}\!\!\cdot F))(a_1,\dots, a_n).$$
 By hypothesis the left side is zero,  so if $r\ne 0$ this shows that $F= (\delta(\iota_{a_0}\!\!\cdot r^{-1} F))(a_1,\dots, a_n)$, a coboundary, while if $r=0$ then it shows that $\iota_{\sigma}\!\cdot F$ is a cocycle.$\Box$

\begin{corollary} Every cocycle is cohomologous to its homogeneous part of degree zero; therefore, $H^n(\g, M)$ is identical with its homogeneous part of degree zero for every $n$. $\Box$
\end{corollary}

When $F$ is a homogeneous $1$-cochain Theorem \ref{zero} asserts that $\iota{\sigma}\cdot\!F =F(\sigma)$ is in $M^{\g}$.

When there are several commuting \emph{adss} elements whose adjoints have eigenvalues in $k$ then the eigenspaces of the adjoint of any one are invariant under the operations of the others. The decomposition of the cohomology can then be further refined into a direct sum of parts which are homogeneous with respect to all the commuting \emph{adss} elements. This multigrading can be represented, as usual, by a linear function from the space $V$ spanned over $k$ by the \emph{adss} elements to $k$, i.e., by an element of $V^*$. The multidegree of an element is then just the element of $V^*$ which it defines.  Every cocycle is then cohomologous to the part which is of homogeneous of degree zero simultaneously with respect to all the gradings induced by the elements, i.e., for which the corresponding element of $V^*$ is zero.  Choosing a basis for $V$, it is convenient to represent the mulltidegree by a ``degree vector'' whose components are the degrees relative to the respective elements of the basis. One then has

\begin{corollary} If $F$ is an $n$-cocycle and $\sigma_1, \dots, \sigma_r$ are commuting \emph{adss} elements then $\iota_{\sigma_r}\iota_{\sigma_{r-1}}\cdots\iota_{\sigma_1}F$ is an $n-r$-cocycle. $\Box$
\end{corollary}

\section{Cohomology of the Borel subalgebra}  As a particulary simple application (in the case where $M = k$) we show how the theorem can be used to compute the cohomology of the the Borel subalgebra $\b$ of  $\mathfrak{sl}(N, k)$, i.e., the subalgebra consisting of all upper triangular matrices.  The Cartan subalgebra $\h$ of diagonal matrices (of trace zero) is then a commutative  subalgebra of  \emph{adss} elements and the $e_{ij}, i < j$ are simultaneous eigenvectors, with their eigenvalues lying in $k$.  Exactly the same result, namely that the cohomology of $\b$ is the exterior algebra on $\h$, holds for all the simple algebras, and can be proved similarly, except that we have not examined for the others the appropriate constraint on the characteristic. For characteristic zero, the result is a special case of a theorem  of Kostant cf \cite{Kostant:Borel}.  

\begin{theorem} Let $\b$ be the Borel subalgebra of $\mathfrak{sl}(N,k)$ and $\h$ be the Cartan subalgebra.  If the characteristic $p$ of $k$ is zero or greater than $N$ then $H^n(\b, k) = (\bigwedge^n \h)^*$ for all $n$.
\end{theorem}
\noindent\textsc{Proof.} Since the coefficient module is $k$ and therefore homogeneous of degree zero, it is sufficient to show that if $F$ is a homogeneous $n$-cocycle of degree zero then it must vanish when 
any argument is an $e_{ij}$ with $i< j$.  As basis for $\h$ take the matrices $h_i = e_{ii}-e_{i+1,i+1}, \, i = 1, \dots, N-1$ and suppose at first that the characteristic is zero.  The degree vectors for the simple roots $$e_{12}, e_{23}, \dots, e_{N-2,N-1}, e_{N-1,N}$$ are then, respectively $$(2,-1, 0,\dots, 0), (-1,2,-1,0,\dots, 0), \dots, (0,\dots, 0, -1,2,-1), (0,\dots,0,-1,2).$$ The others are obtained as sums of these; that for $e_{ij}$ is the sum of those for $e_{i,i+1}, e_{i+1,i+2},\dots,e_{j-1,j}$. Note that the sum of the entries in the degree vector for a simple root $e_{i,i+1}$ vanishes except for $e_{12}$ and $e_{N-1,N}$; in these cases it is $+1$.  Suppose now that $F$ is not zero when evaluated at some $n$-tuple of homogeneous elements.  The homogeneous elements are either elements of $\h$ or multiples of the $e_{ij}$ with $i<j$. There can be no repetitions since $F$ is alternating. If any of the $e_{ij}, i<j$ appear then the sum of their degree vectors must be the zero vector, else $F$ would not be of degree zero. This is impossible if any $e_{1j}$ or $e_{iN}$ appears. By stripping the first and last entry from each degree vector one is then  effectively reduced to the case of dimension $N-2$, proving the proposition in characteristic zero. For characteristic $p > 0$, it is easy to check that no entry in a sum of degree vectors for roots of the form $e_{1j}$ and $e_{iN}$ can exceed $N$ in absolute value, so the argument continues to hold as long as $p > N$. $\Box$


\begin{thebibliography}{10}

\bibitem{HiltonStammbach}
P.J.~Hilton and U.~Stammbach.
\newblock {\em {A course in Homological Algebra }}, {\em Graduate Texts in Mathematics 4}.
\newblock Springer-Verlag, New York/Heidelberg/Berlin, 1971.

\bibitem{Kostant:Borel}
B.~Kostant.
\newblock {Lie algebra cohomology and the generalized Borel-Weil theorem}.
\newblock {\em Ann. of Math.}, 74:329--387, 1961.


\end{thebibliography}
\end{document}